\documentstyle[12pt]{amsart}

\newtheorem{thm}{Theorem}

\newtheorem{lem}{Lemma}

\newtheorem{prop}{Proposition}

\newtheorem{rem}{Remark}

\begin{document}

\title[Generators of a Picard Modular Group]{Generators of a Picard modular group in two complex dimensions}         

\author{Elisha Falbel}
\address{Institut de Math\'ematiques\\Universit\'e Pierre et Marie Curie}
\email{falbel@@math.jussieu.fr}
\author{G\'abor Francsics}
\address{Department of Mathematics\\Michigan State University}
\email{francsics@@math.msu.edu}       
\author{Peter D. Lax}
\address{Courant Institute\\New York University}
\email{lax@@courant.nyu.edu}
\author{John R. Parker}
\address{Department of Mathematical Sciences\\University of Durham}
\email{j.r.parker@@durham.ac.uk}

\date{September 30, 2009}
\keywords{Complex hyperbolic space, Picard modular groups \\The second author is grateful for the hospitality of the Matematical Sciences Research Institute at Berkeley, and the R\'enyi Mathematical Institute, Budapest.}

\begin{abstract}
The goal of the article is to prove that four explicitly given transformations, two Heisenberg translations, a rotation and an involution generate the Picard modular group with Gaussian integers acting on the two dimensional complex hyperbolic space. The result answers positively a question raised by A. Kleinschmidt and D. Persson.  
\end{abstract}
\maketitle

\section{Introduction}

Our main goal in this article is to give a simple, self-contained proof that four explicitly given transformations, two Heisenberg translations, a rotation and an involution generate the two dimensional Picard modular group with Gaussian integers acting on the two dimensional complex hyperbolic space. The result answers positively a question raised by A. Kleinschmidt and D. Persson \cite{KP} and it is used in the work \cite{BKPP} on string compactifications. The method used in the paper gives a simple algorithm to decompose any transformation in the Picard group as a product of the generators.

The complex hyperbolic space ${\bf CH}^{2}$ is the rank one Hermitian symmetric space of noncompact type, 
$SU(2,1)/S(U(2)\times U(1))$. A standard model of the complex hyperbolic space is the complex unit ball of ${\bf C}^{2}=\{z\in{\bf C}^{2};\ |z|<1\}$ equipped with the Bergman metric $g=\sum_{j,k=1}^{2}g_{j,k}(z)dz_{j}\otimes d\bar{z}_{k}$, where $g_{j,k}=
\hbox{const}\cdot\partial_{j}\overline{\partial}_{k}\log(1-|z|^{2})$. This model is the bounded realization of the Hermitian symmetric space ${\bf CH}^{2}$. 
We shall use the unbounded hyperquadric model of the complex hyperbolic space, that is, $D^{2}=\{z\in {\bf C}^{2};\ \Re e z_{2}>\frac{1}{2}|z_{1}|^{2}\}$. 

The holomorphic automorphism group of ${\bf CH}^{2}$, $\hbox{\bf Aut}({\bf CH}^{2})$, consists of 
rational functions $G=(g_{1}, g_{2}):D^{2}\mapsto D^{2}$,
\begin{equation}
g_{j}(z)=\frac{g_{j+1,1}+g_{j+1,2}z_{1}+g_{j+1,3}z_{2}}{g_{1,1}+g_{1,2}z_{1}+g_{1,3}z_{2}},
\label{fraclin}
\end{equation}
$j=1, 2.$ These automorphisms act linearly in homogeneous coordinates $(\zeta_{0}, \zeta_{1}, \zeta_{2})$, $z_{j}=\frac{\zeta_{j}}{\zeta_{0}}$, $j=1, 2$. The corresponding matrix $G=(g_{jk})_{j,k=1}^{3}$ satisfies the condition 
\begin{equation}
G^{*}CG=C,
\label{e25}
\end{equation}

where 
\begin{equation*}
C\equiv\left(\begin{array}{ccc}
0 & 0& -1 \\
0 & 1 & 0 \\
-1 & 0 & 0
\end{array}\right).
\end{equation*}
The determinant of the matrix $G$ is normalized to be equal to $1$. The matrix $C$ is the matrix of the quadratic form of a defining function of $D^{2}$ written in homogeneous coordinates. More precisely,  
$$D^{2}=\{[\zeta]\in {\bf CP}^{2};\ r(\zeta)\equiv \langle C\zeta, \zeta \rangle =-\zeta_{2}\bar{\zeta}_{0}-\zeta_{0}\bar{\zeta}_{2}+|\zeta_{1}|^{2}<0\}.$$

The Picard modular groups are 
$$SU(2,1;{\cal O}_{d}),$$
 where ${\cal O}_{d}$ is the ring of algebraic integers of the imaginary quadratic extension ${\bf Q}(i\sqrt{d})$ for any positive squarefree integer $d$ (see Holzapfel \cite{H1}). The elements of the ring ${\cal O}_{d}$ can be described easily (Hardy-Wright \cite{HW}):
$${\cal O}_{d}=\left\{
\begin{array}{ll}
{\bf Z}[i\sqrt{d}] & \mbox{ if } d\equiv 1,2\ (\mbox{mod }4)\\
{\bf Z}[\frac{1+i\sqrt{d}}{2}] & \mbox{ if } d\equiv 3\ \ (\mbox{mod }4).
\end{array}
\right. $$ 
It is well known that the ring ${\cal O}_{d}$ is Euclidean for positive square free integer $d$ if and only if $d=1, 2, 3, 7, 11$, \cite{ST}, page 92. The Picard modular groups  $SU(2,1;{\cal O}_{d})$ are discrete holomorphic automorphism subgroups of ${\bf CH}^{2}$.

Geometric and spectral properties of discrete lattices acting on complex hyperbolic spaces attracted a lot of attention recently. 
See for example, the work of Goldman and Parker \cite{GP}, Francsics and Lax \cite{FL1}, \cite{FL2}, \cite{FL3}, Falbel and Parker \cite{FP}, Deraux, Falbel and Paupert \cite{DFP}, Schwartz \cite{Sch}, or the survey paper of Parker \cite{P} on the geometric properties of complex hyperbolic lattices. Spectral properties of the automorphic complex hyperbolic Laplace-Beltrami operator were investigated by Epstein, Melrose and Mendoza \cite{EMM}, Reznikov \cite{R}, Lindenstrauss and Venkatesh \cite{LV}. 
Despite the remarkable progress, several important algebraic, geometric and analytic problems are still open. To name a few, we mention the existence of nonarithmetic lattices, or the existence of embedded eigenvalues in the continuous spectrum of the automorphic Laplace-Beltrami operator. A general construction of a fundamental domain was obtained in \cite{GR} for Lie groups. However the exact algebraic and geometric structure is known explicitly only for very few lattices and fundamental domains in higher dimensions. This is in sharp contrast to the case of the real hyperbolic spaces ${\bf H}^{2}$, and ${\bf H}^{3}$. 
Since the influential work of Mostow \cite{M} it is well known that discrete holomorphic automorphism subgroups acting on the complex hyperbolic space ${\bf CH}^{n}$ are particularly hard to analyze.

\section{Preliminaries}

Three important classes of holomorphic automorphisms are Heisenberg translations, dilations, and rotations. 

The Heisenberg translation by 
$a\in\partial D^{2}$, $N_{a}\in
\hbox{\bf Aut}({\bf CH}^{2})$ is defined as 
$$N_{a}(z_{1},z_{2})=(z_{1}+a_{1},z_{2}+a_{2}+z_{1}\bar{a}_{1}).$$ 
If we write $a=(a_{1},a_{2})=(\gamma, \frac{1}{2}|\gamma|^{2}+ir)$ with $\gamma\in{\bf C}$, $r\in{\bf R}$ then the Heisenberg translation is given by 
$N_{(\gamma,\frac{1}{2}|\gamma|^{2}+ir)}(z_{1},z_{2})=
(z_{1}+\gamma,z_{2}+\frac{1}{2}|\gamma|^{2}+ir+z_{1}\bar{\gamma})$. The corresponding matrix representation is
\begin{equation*}
N_{a}\equiv\left(\begin{array}{ccc}
1 & 0 & 0 \\
a_{1} & 1 & 0 \\
a_{2} & \bar{a}_{1} & 1
\end{array}\right)=
\left(\begin{array}{ccc}
1 & 0 & 0 \\
\gamma & 1 & 0 \\
\frac{1}{2}|\gamma|^{2}+ir & \bar{\gamma} & 1
\end{array}\right).
\end{equation*}
Product of the Heisenberg translations $N_{a}$, $N_{b}$ is the Heisenberg translation 
\begin{equation}
N_{a}\circ N_{b}=N_{(a_{1}+b_{1}, a_{2}+b_{2}+\bar{a}_{1}b_{1})}
\label{H1}
\end{equation}
with the parameter  $(a_{1}+b_{1}, a_{2}+b_{2}+\bar{a}_{1}b_{1})\in\partial D^{2}$. Moreover the inverse of $N_{a}$ is the Heisenberg translation 
\begin{equation}
N_{a}^{-1}=N_{(-a_{1}, -a_{2}+|a_{1}|^{2})}.
\label{H2}
\end{equation}

The holomorphic automorphism of $D^{2}$, 
$$A_{\delta}(z)=(\delta z_{1},\delta^{2}z_{2})$$ 
is called dilation with parameter $\delta>0$. Its matrix representation is 
\begin{equation*}
A\equiv\left(\begin{array}{ccc}
\frac{1}{\delta} & 0 & 0 \\
0 & 1 & 0 \\
0 & 0 & \delta
\end{array}\right).
\end{equation*}

Rotation in the first variable by $e^{i\varphi}$, 
$$M_{e^{i\varphi}}(z_{1},z_{2})=(e^{i\varphi}z_{1},z_{2})$$
 is a holomorphic automorphism of $D^{2}$ with $\varphi\in{\bf R}$. There are three matrices  
\begin{equation*}
M\equiv\left(\begin{array}{ccc}
\beta & 0 & 0 \\
0 & \beta^{-2} & 0 \\
0 & 0 & \beta
\end{array}\right),
\end{equation*}
 $\beta=e^{-i\varphi/3+2\pi ik/3}$, $k=0,1,2$ corresponding to the same rotation. 

The holomorphic involution 
$$J(z_{1},z_{2})=(z_{1}/z_{2}, 1/z_{2})$$ 
will also play significant role. A matrix representation of $J$ is 
\begin{equation*}
J\equiv\left(\begin{array}{ccc}
0 & 0 & -1 \\
0 & -1 & 0 \\
-1 & 0 & 0
\end{array}\right).
\end{equation*}
 Notice that $J^{2}=I$, $J$ maps $\infty$ into $(0,0)$.

Let $z$ be a boundary point of $D^{2}$, i.e. $z\in\partial D^{2}\cup\{\infty\}$. The stabilizer subgroup (isotropy subgroup)  $\Gamma_{z}$ of $z$ contains all the holomorphic automorphisms that leave $z$ fixed, that is $\Gamma_{z}\equiv\{g\in\mbox{SU}(2,1);\ g(z)=z\}$. The stabilizer subgroup of $\infty$ consists of lower triangular matrices, that is 
$${\cal P}\equiv\Gamma_{\infty}=\{P\in\mbox{SL}(3,{\bf C}); P^{*}CP=C,\  p_{12}=p_{13}=p_{23}=0\}.$$
The Langlands decomposition of the stabilizer subgroup of $\infty$, ${\cal P}\equiv\Gamma_{\infty}$ will play important role in our method. 
Any element of the stabilizer subgroup $P\in{\cal P}$ can be decomposed as a product of a Heisenberg translation, dilation, and a rotation:
\begin{eqnarray*}
P=\left(\begin{array}{ccc}
p_{11} & 0 & 0 \\
p_{21} & p_{22} & 0 \\
p_{31} & p_{32} & p_{33}
\end{array}\right)=
NAM=
\left(\begin{array}{ccc}
\frac{\beta}{\delta} & 0 & 0 \\
 \frac{\beta\gamma}{\delta}& \beta^{-2} & 0 \\
\frac{\beta}{\delta}(\frac{1}{2}|\gamma|^{2}+ir) & \bar{\gamma}\beta^{-2} & \beta\delta
\end{array}\right).
\end{eqnarray*}
The parameters $\beta, \gamma\in{\bf C}$, $\delta, r\in{\bf R}$ satisfy the restrictions $|\beta|=1$, $\delta>0$. 

We recall from \cite{FL1}, \cite{FL2} that 
the Langlands decomposition can also be used to parametrize a holomorphic automorphism $G=(g_{jk})_{j,k=1}^{3}$ not in the stabilizer subgroup of the ideal point, $\infty$. Let $N_{G(\infty)}$ be the Heisenberg translation mapping $(0,0)$ into $G(\infty)$.  Then the transformation $P\equiv JN_{G(\infty)}^{-1}G$ belongs to the stabilizer subgroup of $\infty$, so 
\begin{equation}
G=N_{G(\infty)}JP=N_{G(\infty)}JNAM.
\label{decomp}
\end{equation}
  
The transformations $N$ and $P$ in the decomposition of 
$G$ are not necessarily in the Picard modular group $\Gamma\equiv SU(2,1; {\cal O}_{d})$, even if $G\in\Gamma$.  The entries of $N$, $P$ are not necessarily integers in the ring ${\cal O}_{d}$. However 
\begin{equation}
g_{1j}=-p_{3j},\ j=1,2,3, 
\label{e1015}
\end{equation}
and so $p_{3j}$, $j=1,2,3$ are integers in the ring ${\cal O}_{d}$. 


\section{Statement of the results}
It is well known that the modular group $PSL(2,{\bf Z})\equiv SL(2,{\bf Z})/\{\pm I\}$ is generated by the transformations $z\mapsto z+1$ and $z\mapsto -\frac{1}{z}$. A. Kleinschmidt and D. Persson \cite{KP} raised the question of an analogous statement for the Picard modular group $SU(2,1;{\bf Z}[i])$, namely if there is a simple description of $SU(2,1;{\bf Z}[i])$ in terms of generators. Our main result is to give an elementary proof that the four transformations, $N_{(0,1)}$, $N_{(1+i,1)}$, $M_{i}$ and $J$ are sufficient to generate $SU(2,1;{\bf Z}[i])$. An application of this description to instanton corrections in string theory can be found in \cite{BKPP}. 

\begin{thm}
The Picard modular group $SU(2,1; {\bf Z}[i])$ is generated by the Heisenberg translations 
\begin{equation}
N_{(0,i)}\equiv\left(\begin{array}{ccc}
1 & 0 & 0 \\
0 & 1 & 0 \\
i & 0 & 1
\end{array}\right),\ \ 
N_{(1+i,1)}
\equiv
\left(\begin{array}{ccc}
1 & 0 & 0 \\
1+i & 1 & 0 \\
1 & 1-i & 1
\end{array}\right),
\label{gen12}
\end{equation}
the rotation by $i$ in the first coordinate
\begin{equation}
M_{i}\equiv\left(\begin{array}{ccc}
i & 0 & 0 \\
0 & -1 & 0 \\
0 & 0 & i
\end{array}\right),
\label{gen3}
\end{equation}
and the involution 
\begin{equation}
J\equiv\left(\begin{array}{ccc}
0 & 0 & -1 \\
0 & -1 & 0 \\
-1 & 0 & 0
\end{array}\right).
\label{gen4}
\end{equation}
\label{t5}
\end{thm}

\begin{rem}
The method used in the paper is constructive. It gives an algorithm to decompose any transformation in $SU(2,1; {\bf Z}[i])$ as a product of the four generators in  (\ref{gen12}), (\ref{gen3}),  and (\ref{gen4}). The main ingredients are the decomposition (\ref{decomp}) and the Euclidean algorithm. 
\label{r5}
\end{rem}

\begin{rem} 
It would be interesting to know if our method can be extended to the other Euclidean rings ${\cal O}_{d}$, more precisely, to the Picard modular groups $SU(2,1;{\cal O}_{d})$, $d=2,3,7,11$.  Little is known about the geometric and algebraic properties, e.g., explicit fundamental domains, generators, presentations of these Picard modular groups, except in the case $d=3$ (see  Falbel-Parker \cite{FP}). We also mention a related result on quaternionic hyperbolic space: Woodward \cite{W} obtained generators for the Picard modular group $PU(2,1; {\cal H})$ where $\cal H$ is the ring of Hurwitz integral quaternions. 
\label{r10}
\end{rem}

In \cite{FL1}, \cite{FL2} an explicit fundamental domain was constructed for the Picard modular group $SU(2,1; {\bf Z}[i])$. The symmetry of this fundamental domain under a nonholomorphic isometry was used to obtain information about the embedded eigenvalues (Maass cusp forms) of 
the automorphic Laplace-Beltrami operator of the group $SU(2,1; {\bf Z}[i])$ in \cite{FL3}. 
In  \cite{FFP} the authors construct a different fundamental domain, determine the geometric, combinatorial structure of this fundamental domain, and obtain a presentation of the group. See also the work \cite{Y}. 

\begin{rem}
We mention that very little is known about the fundamental domains and combinatorial structure of the higher dimensional Picard modular groups $SU(n,1;{\cal O}_{d})$.  More generally, it is a major challenge in complex hyperbolic geometry to understand the geometric structure of discrete lattices in higher dimensions. 
\label{r15}
\end{rem}

%

\section{Proof}

We start by characterizing the stabilizer subgroup of infinity, $\cal P$, and describing the arithmetic properties of the entries of a transformation in $\cal P$. 
\begin{lem}
Let $G=(g_{jk})\in SU(2,1)$. Then $G\in{\cal P}$ if and only if $g_{13}=0$. 
Moreover, if $d=1$ then $P\in {\cal P}\equiv\Gamma_{\infty}(2,1;{\bf Z}[i])$ if and only if the parameters in the Langlands decomposition of $P$ satisfy the conditions 
\begin{equation}
\delta=1,\ \  \beta=1, i, -1, -i,\ \  r\in{\bf Z}, \ \  \gamma\in{\bf Z}[i],\ \  |\gamma|^{2}\in 2{\bf Z}.
\label{e30}
\end{equation} 

\label{l5}
\end{lem}

\noindent{\bf Proof of Lemma \ref{l5}}. It is well known that if $G\in{\cal P}$ then $G$ is lower triangular, so it is enough to prove the converse. Computing and comparing the entries in the lower right corner of (\ref{e25}) we obtain 
$$-\bar{g}_{33}g_{13}+|g_{23}|^{2}-\bar{g}_{13}g_{33}=0.$$
So $g_{13}=0$ implies that $g_{23}=0$. Similarly, comparing the entries in the third row, second column give the equation 
$$-\bar{g}_{33}g_{12}+\bar{g}_{23}g_{22}-\bar{g}_{13}g_{32}=0.$$
Thus $\bar{g}_{33}g_{12}=0$. However, $g_{13}=g_{23}=0$ and $\det G=1$ exclude that $g_{33}=0$. Therefore $g_{12}=0$ and $G$ is lower triangular. 
Let $P\in\Gamma_{\infty}(2,1; {\bf Z}[i])$. Since $p_{11}=\beta/\delta$, $p_{33}=\beta\delta$ are nonzero Gaussian integers with $|\beta|=1$, $\delta>0$ it follows immediately that $\delta=1$, and $\beta=\pm 1,\pm i$. Moreover $p_{31}/\beta=\frac{1}{2}|\gamma|^{2}+ir$ and $p_{21}/\beta=\gamma$ are also Gaussian integers. This proves the second part of Lemma \ref{l5}.

\begin{prop}
The stabilizer subgroup of infinity, $\cal P$, in the Picard modular group $SU(2,1; {\bf Z}[i])$ is generated by the Heisenberg translations $N_{(0,i)}$, $N_{(1+i,1)}$ and the rotation $M_{i}$. 
\label{p5}
\end{prop}

\noindent{\bf Proof of Proposition \ref{p5}}. 
Let $P$ be an element in the stabilizer subgroup $\cal P$. We know that $P$ is lower triangular. According to Lemma \ref{l5} there is no dilation component in its Langlands decomposition, that is  
\begin{eqnarray*}
P=NM=\left(\begin{array}{ccc}
1 & 0 & 0 \\
\gamma & 1 & 0 \\
\frac{1}{2}|\gamma|^{2}+ir & \bar{\gamma} & 1
\end{array}\right)
\left(\begin{array}{ccc}
\beta & 0 & 0 \\
0 & \beta^{-2} & 0 \\
0 & 0 & \beta
\end{array}\right).
\end{eqnarray*}
Since $\beta^{4}=1$, the rotation in $P$ is $M_{i}$, $M_{-1}=M_{i}^{2}$, $M_{-i}=M_{i}^{3}$, or $I=M_{i}^{4}$. Therefore the rotation component of $P$ in the Langlands decomposition is generated by $M_{i}$. 

According to (\ref{H1}) the Heisenberg translation part of $P$ splits as
\begin{equation}
N_{(\gamma, \frac{1}{2}|\gamma|^{2}+ir)}=N_{(0,ri)}\circ N_{(\gamma, \frac{1}{2}|\gamma|^{2})}.
\label{e32}
\end{equation}
Here $N_{(0,ri)}$ can be written as 
\begin{equation}
N_{(0,ri)}=N_{(0,i)}^{r},\ \ r\in{\bf Z}
\label{e35}
\end{equation}
observing that the inverse of $N_{(0,i)}$ is $N_{(0,-i)}$. 

The next step is to decompose $N_{(\gamma, \frac{1}{2}|\gamma|^{2})}$ as a product of Heisenberg translations in the directions $(1+i, 1)$, $(-1+i,1)$ and $(0,i)$. Let $\gamma=m+in\in{\bf Z}[i]$. Then $|\gamma|^{2}=m^{2}+n^{2}\in2{\bf Z}$ according to (\ref{e30}). This means that $m$ and $n$ have the same parity. Therefore we can write $\gamma$ as 
$$\gamma=k(1+i)+l(-1+i)$$ 
with $k\equiv\frac{m+n}{2}\in{\bf Z}$, $l\equiv\frac{n-m}{2}\in{\bf Z}$ and $|\gamma|^{2}=2(k^{2}+l^{2})$. 
Since $(0,-2kli)\in \partial D^{2}$, $((k(1+i),k^{2})\in \partial D^{2}$, and $(l(-1+i),l^{2})\in \partial D^{2}$, it follows from (\ref{H1}), (\ref{e35}) that 
\begin{eqnarray}
N_{(\gamma, \frac{1}{2}|\gamma|^{2})} & = & N_{(0,-2kli)}\circ N_{(\gamma, \frac{1}{2}|\gamma|^{2}+2kli)}\nonumber\\
 \ & = & N_{(0,-2kli)}\circ N_{(k(1+i)+l(-1+i), k^{2}+l^{2}+2kli)}\nonumber\\
 \ & = & N_{(0,-2kli)}\circ N_{(k(1+i),k^{2})}\circ N_{(l(-1+i),l^{2})}\nonumber\\
  \ & = & N_{(0,i)}^{-2kl}\circ N_{(k(1+i),k^{2})}\circ N_{(l(-1+i),l^{2})}.
\label{e40}
\end{eqnarray}
Here the inverse of $N_{(1+i,1)}$ is $N_{(-1-i,1)}$ and an easy induction argument shows that
\begin{equation}
N_{(k(1+i), k^{2})}=N_{(1+i,1)}^{k}
\label{e45}
\end{equation}
for any $k\in{\bf Z}$. Moreover the third factor in (\ref{e40}) can be written as
\begin{equation}
N_{(l(-1+i), l^{2})}=M_{i}\circ N_{(l(1+i),l^{2})}\circ M_{i}^{-1}=
M_{i}\circ N_{(1+i,1)}^{l}\circ M_{i}^{-1}
\label{e50}
\end{equation}
using (\ref{e45}). Combining (\ref{e32}), (\ref{e35}), (\ref{e40}), (\ref{e45}) and (\ref{e50}) we obtain that 
\begin{eqnarray*}
N_{(\gamma, \frac{1}{2}|\gamma|^{2}+ir)} & = & N_{(0,i)}^{r}\circ N_{(0,i)}^{-2kl}\circ 
N_{(1+i,1)}^{k}\circ M_{i}\circ N_{(1+i,1)}^{l}\circ M_{i}^{-1}\\
\ & = & N_{(0,i)}^{r-2kl}\circ 
N_{(1+i,1)}^{k}\circ M_{i}\circ N_{(1+i,1)}^{l}\circ M_{i}^{-1}.
\end{eqnarray*}
This completes the proof of Proposition \ref{p5}. 

\noindent{\bf Proof of Theorem \ref{t5}.} Let $G=(g_{jk})_{j,k=1}^{3}$ be an element of the group $SU(2,1;{\bf Z}[i])$. We may assume that $G$ does not belong to the stabilizer subgroup of infinity, $\cal P$. Then $g_{13}\not=0$ and $G$ maps infinity to $(g_{23}/g_{13}, g_{33}/g_{13})$. Since $G(\infty)$ is in 
$\partial D^{2}$
\begin{equation}
\Re e \frac{g_{33}}{g_{13}}=\frac{1}{2}\big|\frac{g_{23}}{g_{13}}\big|^{2}.
\label{e55}
\end{equation} 
Consider the Heisenberg translation $N_{G(\infty)}$ that maps $(0,0)$ to $G(\infty)$. Note that the translation $N_{G(\infty)}$ is not necessarily in the Picard modular group $SU(2,1;{\bf Z}[i])$ except if $|g_{13}|=1$. However, we know from (\ref{decomp}) that 
$$JN_{G(\infty)}^{-1}G=P.$$
We will successively approximate $N_{G(\infty)}^{-1}$ by Heisenberg translations in the Picard group to decrease the value $|g_{13}|^{2}\in{\bf Z}$ until it becomes $0$. Then $G$ belongs to the stabilizer subgroup $\cal P$ according to Lemma \ref{l5} and can be expressed as a product of the generators (\ref{gen12}) and 
(\ref{gen3}) according to Proposition \ref{p5}. The approximation step uses the fact that the ring ${\cal O}_{1}\equiv{\bf Z}[i]$ is Euclidean. 

Write 
\begin{equation}
-\frac{g_{23}}{g_{13}}=x(1+i)+y(-1+i)=(x-y)+(x+y)i
\label{e60}
\end{equation}
with real numbers $x,y\in{\bf R}$. Select integers $m,n\in{\bf Z}$ such that 
\begin{equation}
|x-m|\leq 1/2,\ \  |y-n|\leq 1/2,
\label{e62}
\end{equation} 
i.e., the nearest integers in $\bf R$. Let 
\begin{equation}
\gamma=m-n+i(m+n), 
\label{e63}
\end{equation}
and select an integer $k\in{\bf Z}$ such that 
\begin{equation}
\big|k+\Im m (\bar{\gamma}\frac{g_{23}}{g_{13}})+\Im m \frac{g_{33}}{g_{13}} \big|\leq\frac{1}{2} .
\label{e65}
\end{equation}
Then we approximate $N_{G(\infty)}^{-1}$ by the translation 
$$N_{(\gamma,\frac{1}{2}|\gamma|^{2}+ik)}\equiv N_{(m-n+i(m+n), m^{2}+n^{2}+ik)}.$$
Notice that $N_{(\gamma,\frac{1}{2}|\gamma|^{2}+ik)}$ is in the Picard modular group because 
\begin{equation}
|\gamma|^{2}=(m-n)^{2}+(m+n)^{2}=2m^{2}+2n^{2}\in 2{\bf Z}.
\label{e70}
\end{equation} 
Next we calculate the entry in the upper right corner of the product 
\begin{eqnarray}
G_{1} & \equiv & JN_{(\gamma,\frac{1}{2}|\gamma|^{2}+ik)}G
\label{e72}\\
\  & = & \left(\begin{array}{ccc}
0 & 0 & -1 \\
0 & -1 & 0 \\
-1 & 0 & 0
\end{array}\right)
\left(\begin{array}{ccc}
1 & 0 & 0 \\
\gamma & 1 & 0 \\
\frac{1}{2}|\gamma|^{2}+ik & \bar{\gamma} & 1
\end{array}\right)G\nonumber\\
\ & = & \left(\begin{array}{ccc}
-\frac{1}{2}|\gamma|^{2}-ik & -\bar{\gamma} & -1 \\
-\gamma & -1 & 0 \\
-1 & 0 & 0
\end{array}\right)G.\nonumber
\end{eqnarray}

So $g_{13}^{(1)}$, the entry in the upper right corner of $G_{1}=(g_{jk}^{(1)})$, is equal to 
\begin{eqnarray}
g_{13}^{(1)} & \equiv &  -(\frac{1}{2}|\gamma|^{2}+ik)g_{13}-\bar{\gamma}g_{23}-g_{33}\nonumber\\
\ & = &  -g_{13}\big(\frac{1}{2}|\gamma|^{2}+ik+\bar{\gamma}\frac{g_{23}}{g_{13}}+\frac{g_{33}}{g_{13}}\big)\nonumber\\
\ & = & -g_{13}\big[\frac{1}{2}|\gamma|^{2}+\Re e (\bar{\gamma}\frac{g_{23}}{g_{13}})+\Re e\frac{g_{33}}{g_{13}}\big]+\nonumber\\
\ & \ & -ig_{13}\big[k+\Im m (\bar{\gamma}\frac{g_{23}}{g_{13}})+\Im m \frac{g_{33}}{g_{13}}\big]\nonumber\\
\ & \equiv & -g_{13}(I_{1}+iI_{2}).
\label{e75}
\end{eqnarray}

Using (\ref{e55}), (\ref{e60}), the definition of $\gamma$ in (\ref{e63}), and (\ref{e70}) we can simplify $I_{1}$: 
%
%
\begin{eqnarray*}
I_1 & = & \frac{1}{2}|\gamma|^2
+\Re\left(\overline{\gamma}\,\frac{g_{23}}{g_{13}}\right)
+\Re\frac{g_{33}}{g_{13}} \\
& = & \frac{1}{2}|\gamma|^2
+\Re\left(\overline{\gamma}\,\frac{g_{23}}{g_{13}}\right)
+\frac{1}{2}\left|\frac{g_{23}}{g_{13}}\right|^2 \\
& = & \frac{1}{2}\left|\gamma+\frac{g_{23}}{g_{13}}\right|^2 \\
& = & \frac{1}{2}\Bigl|(m-n)+(m+n)i-(x-y)-(x+y)i\Bigr|^2 \\
& = & (x-m)^2+(y-n)^2.
\end{eqnarray*}

Then (\ref{e62}) gives the upper bound 
\begin{equation}
|I_{1}|\leq \left(\frac{1}{2}\right)^{2}+\left(\frac{1}{2}\right)^{2}=\frac{1}{2}.
\label{e80}
\end{equation}
Selection of $k$ in (\ref{e65}) gives the inequality $|I_{2}|\leq 1/2$ for the second term in (\ref{e75}). 
Therefore we can estimate $g_{13}^{(1)}$ by combining (\ref{e75}) with (\ref{e65}) and (\ref{e80}):

$$|g_{13}^{(1)}|^{2}=|g_{13}|^{2}|I_{1}+iI_{2}|^{2}=|g_{13}|^{2}(I_{1}^{2}+I_{2}^{2})\leq |g_{13}|^{2}\left[\left(\frac{1}{2}\right)^{2}+\left(\frac{1}{2}\right)^{2}\right]=\frac{1}{2}|g_{13}|^{2}.$$

Repeating this approximation procedure finitely many times we reduce the matrix of the transformation $G$ to the matrix of a transformation $G_{n}$ with $g_{13}^{(n)}=0$. However, according to Lemma \ref{l5}, this condition implies that the $G_{n}$ belongs to the stabilizer subgroup of infinity $\cal P$. Proposition \ref{p5} guarantees that $G_{n}$ is generated by (\ref{gen12}), and (\ref{gen3}).  Since the approximation procedure (\ref{e72}) uses the transformation $J$ and transformations in $\cal P$, Proposition \ref{p5} implies that $G$ is generated by (\ref{gen12}), (\ref{gen3}), and (\ref{gen4}). This completes the proof of Theorem \ref{t5}.

\end{document}